\documentclass[hidelinks,onefignum,onetabnum]{siamart251216}

\usepackage{amsfonts,amssymb,array,bm,booktabs}
\usepackage{algorithmic}
\usepackage{placeins}
\usepackage{needspace}
\usepackage{mathtools}

\newcommand{\standardfigurewidth}{0.90\linewidth}

\newcommand{\R}{\mathbb{R}}
\DeclareMathOperator{\curl}{curl}
\newsiamremark{remark}{Remark}

\headers{de Rham complex structure in Eig-PIELM}{H. Kang et al.}

\title{Embedding a de Rham Complex Structure into Eig-PIELM for Maxwell Eigenvalue Problems \thanks{
\funding{This work is supported by NSFC (Grant No. 12571399 and 12101511) and Fundamental Research Funds for the Central Universities (Grant No. JBK22YJ37).}}}

\author{Hao Kang, Liang Li\thanks{School of Mathematical Sciences, University of Electronic Science and Technology of China (UESTC), Chengdu, China. (\email{kanghao2022@gmail.com, plum.liliang@gmail.com}).}
\and Xing Li\thanks{School of Electronic Science and Engineering, UESTC (\email{xingli@uestc.edu.cn})}
}

\ifpdf
\hypersetup{
  pdftitle={Embedding a de Rham Complex Structure into Eig-PIELM for Maxwell Eigenvalue Problems},
  pdfauthor={Hao Kang}
}
\fi

\begin{document}

\maketitle

\begin{abstract}
The perfect electric conductor (PEC) boundary condition alone does not eliminate the gradient nullspace of the Maxwell curl–curl operator. Consequently, a boundary-only Maxwell implementation of Eig-PIELM, an eigenvalue extension of the physics-informed extreme learning machine (PIELM), may retain represented gradient modes in its reduced trial space, producing nonphysical, gradient-contaminated modes. We develop dR-Eig-PIELM\footnote{The source code repository is available at \url{https://github.com/kanghao2022/dR-Eig-PIELM}.}, a finite-dimensional structure-preserving extension that retains the linear-coefficient, boundary-null-space, and one-shot multimode architecture of Eig-PIELM. It combines compatible tensor-product grad--curl trial spaces with a constrained null-space construction to restrict the sampled PEC space to the mass-orthogonal complement of the represented gradient image before the eigensolve. Symmetric weak-form Gram assembly yields a symmetric-definite reduced eigenproblem without an eigenvalue cutoff under the stated sampling and material-positivity assumptions. For rectangles and cuboids, we establish fixed-degree exactness at the PEC vector-field space and derive rank and dimension identities. Numerical experiments show that dR-Eig-PIELM excludes the spurious eigenmodes with significant gradient components observed in the boundary-only baseline. The method exhibits rapid convergence with increasing polynomial degree and competitive work--precision performance in the tested configurations. A covariantly mapped waveguide example illustrates the potential of dR-Eig-PIELM for integration into design optimization workflows through a constrained local geometry update driven by computed cutoff frequencies and their sensitivities. This work provides a mathematical and computational framework for separating represented gradient modes from physical Maxwell eigenmodes within Eig-PIELM and supports its use in waveguide spectral analysis and geometry design.
\end{abstract}



\section{Introduction}\label{sec:introduction}

Maxwell eigenvalue problems determine cavity resonances, waveguide cutoff
spectra, and electromagnetic field patterns
\cite{BarrenecheaEtAl2014Enclosures,HuangEtAl2015NullSpaceFree}. Their main
discretization difficulty is not the PEC boundary condition alone, but the
gradient kernel of the curl--curl operator. Let $\Omega\subset\R^d$,
$d\in\{2,3\}$, be a bounded domain occupied by a lossless isotropic medium
with positive electric permittivity $\epsilon$ and inverse magnetic
permeability $\nu$. Under a PEC boundary
condition, the electric formulation seeks
$\lambda=\omega^2>0$ and a nonzero field
$\bm E\in H_0(\curl;\Omega)$ that formally satisfies
\cite{Hiptmair2002,Monk2003}
\begin{equation}\label{eq:maxwell-eigenproblem}
  \begin{aligned}
    \nabla\times\!\left(\nu\,\nabla\times\bm E\right)
    &=\lambda\epsilon\bm E
    &&\text{in }\Omega,\\
    \nabla\cdot(\epsilon\bm E)&=0
    &&\text{in }\Omega,\\
    \bm n\times\bm E&=\bm0
    &&\text{on }\partial\Omega.
  \end{aligned}
\end{equation}
On the domains considered here, the curl--curl kernel is
$\nabla H_0^1(\Omega)$ \cite{Hiptmair2002}. These
gradient fields already satisfy the PEC tangential-trace condition. The
boundary condition therefore does not separate them from the physical
subspace \cite{HuangEtAl2015NullSpaceFree}. If a discretization misrepresents
this kernel, gradient-like directions can appear as spurious positive
eigenmodes or cause spectral pollution
\cite{BarrenecheaEtAl2014Enclosures,Boffi1999Resonators,Boffi2010,Hiptmair2002}.

A Maxwell discretization must therefore represent the curl kernel as the
image of a discrete gradient operator and separate that image from the
physical eigenspace. Finite element exterior calculus (FEEC) formulates this
requirement through discrete de Rham complexes and commuting projections
\cite{Arnold2006FEEC}. N\'ed\'elec edge elements provide the classical
$H(\curl)$-conforming realization \cite{Hiptmair2002,Nedelec1980}. Exactness
controls the kernel identity at a fixed discretization level. Approximation
and discrete compactness are additionally needed for asymptotic spectral
correctness \cite{ArnoldAwanou2014Cubical,Boffi2010,BoffiEtAl2011DiscreteCompactness,Hiptmair2002,Kikuchi1989}.

Existing methods enforce this requirement either through their trial spaces
or through algebraic treatment of the kernel. Compatible spline and cubical
complexes preserve the de Rham sequence, and Bernstein representations reduce
the cost of high-order differential-form calculations
\cite{BuffaEtAl2011IGAForms,BuffaSangalliVazquez2010,
Kirby2014LowComplexity}. Null-space-free eigensolvers use discrete vector
potentials to keep iterations outside the gradient kernel
\cite{HuangEtAl2015NullSpaceFree}. Discrete involutions and
kernel-compensated mimetic finite differences provide related treatments in
other discretization settings
\cite{ErnGuermond2023Involution,ErnGuermond2025Maxwell,
JinXiaXu2025KernelCompensation}. Although their realizations differ, these
methods all control the relation among the gradient image, the curl kernel,
and the physical eigenspace.

The same structural questions arise in neural and random-feature
discretizations of Maxwell eigenproblems.

Neural PDE solvers optimize residuals, energy functionals, or Rayleigh
quotients \cite{EYu2018DeepRitz,Jin2022PINNEigen,PangEtAl2019fPINN,Raissi2019PINN,Rowan2025Rayleigh}. Related collocation and tensor methods use
adaptive basis enrichment or variational assembly
\cite{WangXie2024MultiEigen,WengMaoShen2026DeepCollocation}. For Maxwell
eigenproblems, FieldTNN uses a trainable tensor representation with a
divergence penalty \cite{Jiang2026FieldTNN}. Other neural discretizations
represent finite-element de Rham spaces or compatible differential forms
\cite{Badia2025CompatibleFEINN,Hu2026NeuralDeRham,Longo2023DeRhamNN}.

Random-feature and extreme-learning methods fix the feature map and determine
only linear output coefficients
\cite{Chen2022RFM,Dwivedi2020PIELM,Huang2006ELM,NelsenStuart2021RandomFeature,
Xiong2025RFMEigen}. Eig-PIELM combines this representation
with a sampled boundary null space and one reduced multimode eigensolve, providing an efficient way to compute multiple eigenpairs \cite{Mishra2026EigPIELM}.

However, a direct Maxwell specialization of Eig-PIELM does not
automatically inherit the compatible structure described above.
Choosing scalar features independently for each vector component need
not produce a compatible grad--curl sequence. Moreover, the sampled
boundary null space enforces the PEC condition but does not separate
the represented gradient image from the physical eigenspace, because
these gradient fields already satisfy the tangential boundary
condition. Thus, boundary enforcement alone does not resolve the
structural difficulty of the Maxwell eigenproblem.

The challenge is to incorporate this structure into the
linear-coefficient architecture of Eig-PIELM. The scalar and vector
trial spaces must be coordinated so that the represented gradient
image lies in the PEC-admissible vector space and, for the spaces
considered here, coincides with its curl kernel. This image must then
be separated in the material-weighted mass inner product. The resulting
constraints should define a reduced physical-space eigenproblem while
retaining the one-shot multimode solution procedure of Eig-PIELM.

To address these issues, we introduce de Rham compatibility into the
trial-space construction of Eig-PIELM while retaining its
linear-coefficient, boundary-null-space, and one-shot multimode
architecture. A compatible grad--curl sequence links the scalar and
vector trial spaces. Within the vector space, we impose the sampled
tangential PEC condition and construct a second null space to remove
the represented gradient image in the mass inner product. The
eigenproblem is thus posed directly on the discrete Maxwell physical
space. Symmetric weak-form Gram assembly yields a symmetric-definite
reduced eigenproblem under the stated sampling and
material-positivity assumptions.

\begin{samepage}
The paper makes three main contributions.
\begin{itemize}
  \item We formulate dR-Eig-PIELM as a finite-dimensional
  structure-preserving extension of Eig-PIELM for PEC Maxwell eigenproblems.
  The constrained null-space basis satisfies
  \[
    \operatorname{range}\mathbf Z_X
    =\ker\mathbf B_{\mathrm{bc}}
    \cap
    \bigl(\operatorname{range}\bm\Gamma\bigr)^{\perp_{\mathbf M}}.
  \]
  This basis realizes the discrete Maxwell physical space before the
  eigensolve and excludes the represented gradient image without an
  eigenvalue cutoff. The construction retains the linear-coefficient and
  one-shot multimode architecture of Eig-PIELM.

  \item For tensor-product polynomial spaces on rectangles and cuboids, we
  establish PEC compatibility and fixed-degree exactness at the vector-field
  space,
  \[
    \ker\!\left(
      \left.\nabla\times\right|_{\mathbb V_{p,d}^{0}}
    \right)
    =\nabla\mathbb S_{p,d}^{0}.
  \]
  We derive the associated rank and dimension identities. Under the stated
  trace-unisolvence, evaluation-injectivity, and material-positivity
  assumptions, the restricted mass and curl matrices are positive definite
  and define a well-posed symmetric generalized eigenproblem.

  \item Numerical experiments compare dR-Eig-PIELM with the
  boundary-only Eig-PIELM baseline, analytic spectra, and independent
  N\'ed\'elec references. They verify the tested structural identities
  and show that the method excludes the spurious eigenmodes with
  significant gradient components observed in the boundary-only
  baseline. The tested configurations exhibit rapid convergence with
  increasing polynomial degree and competitive work--precision
  performance. 
\end{itemize}
\end{samepage}

The remainder of the paper is organized as follows.
Section~\ref{sec:maxwell-compatible-structure} introduces the Maxwell setting,
PIELM, and the Eig-PIELM baseline. Section~\ref{sec:algebraic-algorithm}
presents the construction and assembly of dR-Eig-PIELM and its numerical
realization, while Section~\ref{sec:compatible-analysis} establishes its
structural analysis. Section~\ref{sec:numerical-experiments} presents the
numerical results. Section~\ref{sec:conclusion} concludes the paper, with
implementation and reproducibility details deferred to the appendices.

\section{Mathematical Foundations}
\label{sec:maxwell-compatible-structure}

We first define the continuous and discrete Maxwell physical spaces, then review the physics-informed extreme learning machine (PIELM) and the Eig-PIELM baseline.

\subsection{Maxwell Eigenvalue Problem}

\subsubsection{Continuous setting}

Let $d\in\{2,3\}$, and let $\Omega\subset\R^d$ be a bounded Lipschitz domain. We consider the source-free Maxwell eigenproblem \eqref{eq:maxwell-eigenproblem} with $\lambda=\omega^2$, under the
coefficient assumptions stated in Section 1. 

Let $\gamma_t$ and $\gamma_0$ denote the tangential and scalar trace
operators, respectively, with
$\gamma_t(\bm v)=\bm n\times\bm v$ and
$\gamma_0(q)=q|_{\partial\Omega}$ for smooth fields
\cite{Monk2003,Hiptmair2002}. The associated
homogeneous trace spaces are
\[
\begin{aligned}
	H_0(\curl;\Omega)
	&\coloneqq \ker\gamma_t
	=\{\bm v\in H(\curl;\Omega):\gamma_t(\bm v)=\bm 0\},\\
	H_0^1(\Omega)
	&\coloneqq \ker\gamma_0
	=\{q\in H^1(\Omega):\gamma_0(q)=0\}.
\end{aligned}
\]

For $\bm u,\bm v\in H_0(\curl;\Omega)$, define
\[
	\begin{aligned}
		a_\nu(\bm u,\bm v)
		&=
		\int_\Omega
		\nu\,(\nabla\times\bm u)\cdot(\nabla\times\bm v)
		\,\mathrm d\bm x,
		\qquad
		m_\epsilon(\bm u,\bm v)
		&=
		\int_\Omega
		\epsilon\,\bm u\cdot\bm v
		\,\mathrm d\bm x.
	\end{aligned}
\]
The standard weak Maxwell eigenproblem is \cite{Monk2003,Hiptmair2002}
\begin{equation}\label{eq:continuous-weak-eigenproblem}
	\begin{aligned}
		a_\nu(\bm E,\bm v)
		&=
		\lambda m_\epsilon(\bm E,\bm v)
		&&\text{for all }\bm v\in H_0(\curl;\Omega),\\
		m_\epsilon(\bm E,\bm E)&=1.
	\end{aligned}
\end{equation}

For $q\in H_0^1(\Omega)$, one has
$\nabla q\in H_0(\curl;\Omega)$ and
$\nabla\times\nabla q=0$. These mapping properties are part of the
relative de Rham structure \cite{Arnold2006FEEC,Hiptmair2002}. The relevant
portions of the corresponding complexes are
\begin{subequations}\label{eq:continuous-relative-derham-complexes}
\begin{equation}\label{eq:continuous-relative-derham-2d}
H_0^1(\Omega)
\xrightarrow{\ \nabla\ }
H_0(\curl;\Omega)
\xrightarrow{\ \nabla\times\ }
L_0^2(\Omega),
\quad d=2,
\end{equation}
\begin{equation}\label{eq:continuous-relative-derham-3d}
H_0^1(\Omega)
\xrightarrow{\ \nabla\ }
H_0(\curl;\Omega)
\xrightarrow{\ \nabla\times\ }
H_0(\operatorname{div};\Omega)
\xrightarrow{\ \nabla\cdot\ }
L_0^2(\Omega),
\quad d=3.
\end{equation}
\end{subequations}

Since consecutive differential operators compose to zero, this gives
\[
\nabla H_0^1(\Omega)
\subset
\ker\!\left(
\left.\nabla\times\right|_{H_0(\curl;\Omega)}
\right).
\]
If, in addition, $\Omega$ is contractible and $\partial\Omega$ is connected,
the corresponding relative de Rham sequence is exact at
$H_0(\curl;\Omega)$ \cite{Hiptmair2002,Arnold2006FEEC}, and the inclusion
becomes
\begin{equation}\label{eq:continuous-curl-kernel}
\ker\!\left(
\left.\nabla\times\right|_{H_0(\curl;\Omega)}
\right)
=
\nabla H_0^1(\Omega).
\end{equation}
For every $q\in H_0^1(\Omega)$,
\[
\nabla q\in H_0(\curl;\Omega),
\qquad
\nabla\times\nabla q=0.
\]
Hence, $\nabla H_0^1(\Omega)$ is a subspace whose
fields already satisfy the tangential PEC condition. The boundary condition
alone therefore does not separate the gradient kernel.

For a positive eigenpair, testing
\eqref{eq:continuous-weak-eigenproblem} with $\bm v=\nabla q$ gives
\[
0=a_\nu(\bm E,\nabla q)
=\lambda m_\epsilon(\bm E,\nabla q),
\qquad q\in H_0^1(\Omega).
\]
Thus every positive-eigenvalue eigenfunction belongs to
\begin{equation}\label{eq:continuous-physical-space}
\mathbb X_\epsilon
\coloneqq
H_0(\curl;\Omega)
\cap
\bigl(\nabla H_0^1(\Omega)\bigr)^{\perp_{m_\epsilon}}.
\end{equation}
This $m_\epsilon$-orthogonality imposes the source-free Gauss constraint
$\nabla\cdot(\epsilon\bm E)=0$
\cite{Monk2003,DuanMaZou2021Gauss}.
Under \eqref{eq:continuous-curl-kernel}, the restriction in
\eqref{eq:continuous-physical-space} removes the entire curl-free PEC kernel.

\subsubsection{Discrete setting}

Following the discrete de Rham framework
\cite{Arnold2006FEEC,Hiptmair2002}, let $\mathbb S$, $\mathbb V$, and
$\mathbb Q$ be finite-dimensional
scalar-potential, vector-field, and curl-target spaces, respectively.
Let $\mathcal G$ and $\mathcal C$ denote the discrete gradient and curl
operators. The discrete counterpart of
\eqref{eq:continuous-relative-derham-complexes} starts with the grad--curl
chain
\begin{equation}\label{eq:discrete-grad-curl-chain}
\mathbb S
\xrightarrow{\ \mathcal G\ }
\mathbb V
\xrightarrow{\ \mathcal C\ }
\mathbb Q,
\qquad
\mathcal C\mathcal G=0.
\end{equation}

We next introduce the homogeneous boundary constraints relevant to the
PEC problem. Let $\mathbb S_0\subset\mathbb S$ denote the homogeneous
scalar-potential subspace corresponding to $H_0^1(\Omega)$. Let
$\mathcal B$ denote the discrete tangential-trace operator on
$\mathbb V$, and define the PEC-admissible vector subspace by
\[
\mathbb V_0
\coloneqq
\ker\mathcal B
\subset\mathbb V.
\]

Compatibility with the PEC constraint requires
\begin{equation}\label{eq:discrete-pec-compatibility}
\mathcal G(\mathbb S_0)\subset\mathbb V_0,
\qquad\Longleftrightarrow\qquad
\mathcal B\mathcal G q=0
\quad\text{for every }q\in\mathbb S_0.
\end{equation}
Combining \eqref{eq:discrete-grad-curl-chain} and
\eqref{eq:discrete-pec-compatibility} gives
$\mathcal G(\mathbb S_0)\subset
\ker(\left.\mathcal C\right|_{\mathbb V_0})$.

The constrained grad--curl structure is exact at $\mathbb V_0$ if
\begin{equation}\label{eq:discrete-curl-kernel-exactness}
\ker\!\left(
\left.\mathcal C\right|_{\mathbb V_0}
\right)
=
\mathcal G(\mathbb S_0).
\end{equation}
This is the discrete counterpart of the continuous kernel identity
\eqref{eq:continuous-curl-kernel}.

Let $m$ denote a symmetric discrete mass form on $\mathbb V$ whose restriction
to $\mathbb V_0$ is positive definite. The corresponding discrete physical
space is
\begin{equation}\label{eq:abstract-physical-space}
\mathbb X
\coloneqq
\left\{
\bm v\in\mathbb V_0:
m\bigl(\bm v,\mathcal G q\bigr)=0
\quad
\text{for all }q\in\mathbb S_0
\right\}
=
\mathbb V_0
\cap
\bigl(\mathcal G(\mathbb S_0)\bigr)^{\perp_m}.
\end{equation}
It is the discrete analogue of \eqref{eq:continuous-physical-space}. Since
$\mathcal G(\mathbb S_0)\subset\mathbb V_0$ and $m$ is positive definite on
$\mathbb V_0$, this definition gives the $m$-orthogonal direct sum
\begin{equation}\label{eq:abstract-physical-decomposition}
\mathbb V_0
=
\mathcal G(\mathbb S_0)
\oplus_m
\mathbb X.
\end{equation}
Under \eqref{eq:discrete-curl-kernel-exactness}, the first summand is the
entire discrete curl-free PEC kernel.

\subsection{Physics-Informed Extreme Learning Machine}
\label{sec:pielm}

The extreme learning machine (ELM) uses a single-hidden-layer network whose
hidden parameters are sampled once and then fixed. Only the output weights are
determined from data, usually by a linear least-squares solve
\cite{Huang2006ELM}. PIELM applies this fixed-feature construction to
differential equations \cite{Dwivedi2020PIELM}. For scalar notation, it writes
the trial function as
\begin{equation}\label{eq:pielm-trial-function}
  u_N(\bm x)=\sum_{j=1}^{N_\phi}\beta_j\phi_j(\bm x),
\end{equation}
where the hidden features $\phi_j$ are fixed and the coefficients $\beta_j$
are unknown.

For a linear boundary-value problem $\mathcal L u=f$ in $\Omega$ and
$\mathcal B u=g$ on $\partial\Omega$, PIELM evaluates the differential and
boundary operators at collocation points. The two sets of equations form the
linear system
\[
  \begin{bmatrix}
    \bm A_{\Omega}\\
    \bm A_{\partial\Omega}
  \end{bmatrix}\bm\beta
  =
  \begin{bmatrix}
    \bm f_{\Omega}\\
    \bm g_{\partial\Omega}
  \end{bmatrix},
  \qquad
  (\bm A_{\Omega})_{ij}=\mathcal L\phi_j(\bm x_i),
  \quad
  (\bm A_{\partial\Omega})_{ij}=\mathcal B\phi_j(\bm x_i).
\]
An overdetermined system is solved in the least-squares sense. With fixed features, PIELM reduces the PDE solve to linear algebra without backpropagation. 

\subsection{Eig-PIELM Framework}
\label{sec:eigpielm-framework}

Eig-PIELM extends the random fixed-feature PIELM construction to the linear eigenproblem
\cite{Mishra2026EigPIELM}
\begin{equation}
	\label{eq:eigpielm-linearproblem}
\mathcal{L}u=\lambda u \quad \text{in } \Omega,\quad
\mathcal{B}u=0 \quad \text{on } \partial\Omega,\quad
u\ne 0.
\end{equation}
It uses the PIELM ansatz \eqref{eq:pielm-trial-function}.  Both $\bm\beta$ and
$\lambda$ are unknown, so the forced-problem least-squares solve no longer
applies directly.

Let $\{\bm{x}_i\}_{i=1}^{N_x} \subset \Omega$ be interior collocation points.
Evaluate the fixed features and their images under $\mathcal L$ at these points:
\[
	\bm{\Phi}_{ij} = \phi_j(\bm{x}_i),
	\qquad
	\bm{\Psi}_{ij} = (\mathcal{L}\phi_j)(\bm{x}_i).
\]
By linearity of $\mathcal L$, substituting \eqref{eq:pielm-trial-function} into \eqref{eq:eigpielm-linearproblem} gives the sampled residual
$$
  \bigl[(\bm\Psi-\lambda\bm\Phi)\bm\beta\bigr].
$$
Eig-PIELM then forms the coefficient-space products
\[
	\mathbf A_{\mathrm{Eig}}
	=
	\bm{\Phi}^{\top}\bm{\Psi},
	\qquad
	\mathbf M_{\mathrm{Eig}}
	=
	\bm{\Phi}^{\top}\bm{\Phi}.
\]
Eig-PIELM restricts the trial function by the homogeneous boundary
condition before the eigensolve. For boundary points
$\{\bm x_\ell^b\}_{\ell=1}^{N_b}\subset\partial\Omega$, define the
boundary sampling matrix by
$$(\mathbf B_{\mathrm{bc}})_{\ell j}
=(\mathcal B\phi_j)(\bm x_\ell^b),$$ so the sampled boundary condition is
$\mathbf B_{\mathrm{bc}}\bm\beta=0$. An SVD of this matrix gives an orthonormal null-space basis $\mathbf Z_{\mathrm{bc}}$.
Every coefficient vector satisfying the sampled boundary condition has the form $\bm\beta=\mathbf Z_{\mathrm{bc}}\bm y$. 

We require the sampled interior residual to be orthogonal to
$\operatorname{range}(\bm\Phi\mathbf Z_{\mathrm{bc}})$, the sampled boundary-admissible trial space. This gives the reduced generalized eigenproblem:
\[
	\left(\mathbf{\Phi}\mathbf Z_{\mathrm{bc}}\right)^{\top}
	\left(\mathbf{\Psi}-\lambda\mathbf{\Phi}\right)
	\mathbf Z_{\mathrm{bc}}\bm y
	=
	\left(
	\mathbf Z_{\mathrm{bc}}^{\top}\mathbf A_{\mathrm{Eig}}\mathbf Z_{\mathrm{bc}}
	-
	\lambda
	\mathbf Z_{\mathrm{bc}}^{\top}\mathbf M_{\mathrm{Eig}}\mathbf Z_{\mathrm{bc}}
	\right)\bm y
	=0.
\]

Solving this reduced generalized eigenproblem gives $\bm y$. The eigenfunction
is recovered from $\bm\beta=\mathbf Z_{\mathrm{bc}}\bm y$.

A direct componentwise application of Eig-PIELM to a two-dimensional Maxwell
field assigns the same scalar basis independently to the two components:
\begin{equation}\label{eq:eigpielm-trial-field}
	\bm E
	=
	\begin{pmatrix}
		\displaystyle\sum_{j=1}^{N_\phi}c_j^{(1)}\phi_j\\[2mm]
		\displaystyle\sum_{j=1}^{N_\phi}c_j^{(2)}\phi_j
	\end{pmatrix}.
\end{equation}
Because the two component spaces are chosen independently, this trial space
need not realize the discrete grad--curl chain
\eqref{eq:discrete-grad-curl-chain}.  The boundary projection therefore
identifies only $\ker\mathbf B_{\mathrm{bc}}$; it does not represent the
gradient image $\mathcal G(\mathbb S_0)$ or enforce the physical-space
restriction \eqref{eq:abstract-physical-space}.

\section{dR-Eig-PIELM Method}
\label{sec:algebraic-algorithm}

We now realize the grad--curl chain and mass-orthogonal physical space from
Section~\ref{sec:maxwell-compatible-structure} in Eig-PIELM. The construction
proceeds through compatible trial spaces, PEC restriction, and mass-orthogonal
removal of the represented gradient image. The resulting discrete physical
space defines the reduced Maxwell eigenproblem.

\subsection{Construction of compatible trial spaces}
\label{subsec:compatible-operator-realization}

We construct compatible tensor-product trial spaces by matching their
componentwise polynomial degrees to the degree-lowering action of
differentiation.
Let $\mathbb S_{r,s,t}$ denote the tensor product of univariate polynomial
spaces of degrees at most $r$, $s$, and $t$ in the coordinates $x$, $y$, and
$z$. The three-dimensional scalar, vector, curl, and divergence spaces at
level $p$ are
\begin{equation}\label{eq:three-dimensional-first-form-space}
\begin{aligned}
  \mathbb S_p^{(3)}&=\mathbb S_{p,p,p},\\
  \mathbb V_p^{(3)}&=\mathbb S_{p-1,p,p}\times
    \mathbb S_{p,p-1,p}\times\mathbb S_{p,p,p-1},\\
  \mathbb W_p^{(3)}&=\mathbb S_{p,p-1,p-1}\times
    \mathbb S_{p-1,p,p-1}\times\mathbb S_{p-1,p-1,p},\\
  \mathbb Q_p^{(3)}&=\mathbb S_{p-1,p-1,p-1}.
\end{aligned}
\end{equation}
With the component order $(x,y,z)$, they form the sequence
\begin{equation}\label{eq:graded-three-dimensional-sequence}
  \mathbb S_p^{(3)}
  \xrightarrow{\ \nabla\ }
  \mathbb V_p^{(3)}
  \xrightarrow{\ \nabla\times\ }
  \mathbb W_p^{(3)}
  \xrightarrow{\ \nabla\cdot\ }
  \mathbb Q_p^{(3)}.
\end{equation}
Differentiation reduces the polynomial degree in one coordinate by one and
leaves the other degrees unchanged. The componentwise degree shifts in
\eqref{eq:three-dimensional-first-form-space} therefore place the gradient,
curl, and divergence in the successive tensor-product spaces of
\eqref{eq:graded-three-dimensional-sequence}. The grad--curl part realizes the
scalar and electric-field positions of the abstract chain in
Section~\ref{sec:maxwell-compatible-structure}
\cite{Arnold2006FEEC,Hiptmair2002}.

Let $\mathbb S_{r,s}$ be the tensor product in $x$ and $y$. For fields
independent of $z$, take a scalar $q(x,y)$ and a planar vector
$(v_x(x,y),v_y(x,y),0)$. The $z$-component of its curl gives the
two-dimensional sequence
\begin{equation}\label{eq:graded-grad-curl-sequence}
  \mathbb S_{p,p}
  \xrightarrow{\ \mathcal G_p\ }
  \mathbb S_{p-1,p}\times\mathbb S_{p,p-1}
  \xrightarrow{\ \mathcal C_p\ }
  \mathbb S_{p-1,p-1},
\end{equation}
where $\mathcal G_pq=(\partial_xq,\partial_yq)$ and
$\mathcal C_p(v_x,v_y)=\partial_xv_y-\partial_yv_x$.

All dR-Eig-PIELM numerical experiments use the same global tensor-polynomial
basis, evaluated as normalized shifted-Legendre polynomials on the reference
interval $(0,1)$:
\begin{equation}\label{eq:shifted-legendre-reference-basis}
  \widehat\psi_j(\xi)
  =\sqrt{2j+1}\,P_j(2\xi-1),
  \quad 0\leq j\leq p,
\end{equation}
where $P_j$ is the Legendre polynomial on $[-1,1]$. On a physical interval
$(0,L)$, the implementation evaluates $\widehat\psi_j(t/L)$ and scales its
derivatives by $L^{-1}$. This basis spans the same tensor-polynomial spaces as
Bernstein polynomials.

\subsection{Discrete forms}

Let $\bm\Gamma_p$ and $\mathbf R_p$ denote the gradient and curl
coefficient matrices specified in
Appendix~\ref{app:three-dimensional-coefficient-realization}.
In the two-dimensional realization, let $\{\bm\varphi_j\}_{j=1}^{N_v}$
be a basis of $\mathbb S_{p-1,p}\times\mathbb S_{p,p-1}$ and write
\[
  \bm E(\bm c)=\sum_{j=1}^{N_v}c_j\bm\varphi_j.
\]
Choose volume quadrature nodes $\{\bm x_\ell\}$. Evaluating the vector basis
and its curls at these nodes gives
\[
  \begin{aligned}
    (\bm\Phi)_{(\ell,\alpha),j}
    &=[\bm\varphi_j(\bm x_\ell)]_\alpha,
    &\bm\Phi&\in\R^{N_f\times N_v},\\
    \mathbf C_{(\ell,\beta),j}
    &=[\mathcal C_p(\bm\varphi_j)(\bm x_\ell)]_\beta,
    &\mathbf C&\in\R^{N_c\times N_v}.
  \end{aligned}
\]
Here $\alpha$ and $\beta$ index the field and curl components.
Because $\mathbf C$ samples the curl represented by $\mathbf R_p$, the exact
identity $\mathbf R_p\bm\Gamma_p=0$ also yields
$\mathbf C\bm\Gamma_p=0$.

Let $\mathbf W_\epsilon$ and $\mathbf W_\nu$ contain the positive quadrature
and material weights for permittivity and inverse permeability. The weighted
samples are
\[
  \widehat{\bm\Phi}=\mathbf W_\epsilon^{1/2}\bm\Phi,
  \quad
  \widehat{\mathbf C}=\mathbf W_\nu^{1/2}\mathbf C.
\]
Their Gram products define
\[
  \mathbf M=\widehat{\bm\Phi}^{\top}\widehat{\bm\Phi},
  \quad
  \mathbf K=\widehat{\mathbf C}^{\top}\widehat{\mathbf C}.
\]
If $\bm c,\bm d\in\R^{N_v}$ represent two fields in this vector trial space, then
the discrete mass and curl forms are
\[
  m\bigl(\bm E(\bm c),\bm E(\bm d)\bigr)
  =\bm d^{\top}\mathbf M\bm c,
  \quad
  a\bigl(\bm E(\bm c),\bm E(\bm d)\bigr)
  =\bm d^{\top}\mathbf K\bm c.
\]
The sampled weak residual is
\[\bm d^{\top}(\mathbf K-\lambda\mathbf M)\bm c.\]

For the tensor-product realization, the weighted evaluations are assembled by
tensor Gauss quadrature. The weights include the physical area factor and the
sampled material coefficients. Positive weights make $\mathbf M$ and
$\mathbf K$ symmetric positive semidefinite; definiteness on a specified
subspace also requires injective field or curl evaluation there. The fixed
basis in \eqref{eq:shifted-legendre-reference-basis} preserves the
compatible spaces and their exact derivative image.

\subsection{PEC Boundary Constraint}
\label{subsec:pec-admissible-space}

The PEC condition is imposed after the compatible sequence has been
constructed. Define the homogeneous scalar subspace
\begin{equation}\label{eq:homogeneous-scalar-space}
  \mathbb S_{p,p}^0
  =\mathbb S_{p,p}\cap H_0^1(\Omega),
  \quad N_s=\dim\mathbb S_{p,p}^0,
\end{equation}
and let $\mathbf T_0$ inject its coefficient vectors into the coefficients of
$\mathbb S_{p,p}$. The gradient matrix restricted to this subspace is
\[
  \bm\Gamma=\bm\Gamma_p\mathbf T_0
  \in\R^{N_v\times N_s}.
\]
The corresponding PEC vector space is
\[
  \bigl(\mathbb S_{p-1,p}\times\mathbb S_{p,p-1}\bigr)
  \cap H_0(\curl;\Omega).
\]

Let $\mathbf B_{\mathrm{bc}}\in\R^{N_b\times N_v}$ sample the tangential PEC
trace of the vector basis. On the rectangle it evaluates the $x$-component on
horizontal edges and the $y$-component on vertical edges. Because every
function in $\mathbb S_{p,p}^0$ has zero scalar trace, its gradient has zero
tangential trace. Hence
\[
  \mathbf B_{\mathrm{bc}}\bm\Gamma=0.
\]

As in the generic Eig-PIELM boundary-null-space reduction described in
Section~\ref{sec:eigpielm-framework}, a full SVD of
$\mathbf B_{\mathrm{bc}}$ with rank threshold $\tau_b$ supplies a
Euclidean-orthonormal right-null-space basis $\mathbf N_0$. Then
$\bm c=\mathbf N_0\bm y$ parameterizes the sampled PEC-admissible fields. The
sampled null space equals the tangential-trace-zero subspace of the vector
trial space in \eqref{eq:graded-grad-curl-sequence} when boundary sampling is
injective on its represented trace space.

Euclidean coefficient orthogonality does not imply orthogonality in the
sampled mass inner product: $\mathbf N_0^{\top}\mathbf M\mathbf N_0$ is the
Gram matrix of $\widehat{\bm\Phi}\mathbf N_0$ and may be ill-conditioned.
When this sample matrix has full column rank, optional post-PEC
$\mathbf M$-orthonormalization gives dR-Eig-PIELM-stable. We use the same symbol
$\mathbf N_0$ for the PEC basis with or without this step. In either case,
$\operatorname{range}\mathbf N_0=\ker\mathbf B_{\mathrm{bc}}$; after
orthonormalization, $\mathbf N_0^{\top}\mathbf M\mathbf N_0=\mathbf I$.
Algorithm~\ref{alg:drelm-complete} specifies the implementation.

\subsection{Construction of the discrete physical space}

Assume that the gradient map is injective on
$\mathbb S_{p,p}^0$ and that $\mathbf M$ is positive definite on
$\ker\mathbf B_{\mathrm{bc}}$. The basis representation of
\eqref{eq:abstract-physical-space} is
\[
  \mathbb X
  =\ker\mathbf B_{\mathrm{bc}}
  \cap\bigl(\operatorname{range}\bm\Gamma\bigr)^{\perp_{\mathbf M}}.
\]

Let $\{\chi_k^0\}_{k=1}^{N_s}$ be a basis of
$\mathbb S_{p,p}^0$. For $\bm s\in\R^{N_s}$, write
$q(\bm s)=\sum_ks_k\chi_k^0$. A PEC-admissible field belongs to $\mathbb X$
if and only if
\[
  m\bigl(\bm E(\bm c),\mathcal G_p(q(\bm s))\bigr)
  =\bm c^{\top}\mathbf M\bm\Gamma\bm s=0
  \quad\text{for every }\bm s\in\R^{N_s}.
\]
For the PEC basis $\mathbf N_0$, define the mass coupling to the represented
gradients by
\begin{equation}\label{eq:boundary-gradient-coordinates}
  \bm\Gamma_0
  =\mathbf N_0^{\top}\mathbf M\bm\Gamma.
\end{equation}
Substituting $\bm c=\mathbf N_0\bm y$ into the mass-orthogonality condition
gives
\[
  \bm c^{\top}\mathbf M\bm\Gamma\bm s
  =\bm y^{\top}\bm\Gamma_0\bm s
  \quad\text{for every }\bm s\in\R^{N_s}.
\]
Because $\bm s$ is arbitrary, the condition is equivalent to
\[
  \bm\Gamma_0^{\top}\bm y=0,
  \quad \bm y\in\ker\bm\Gamma_0^{\top}.
\]
Same to the boundary constraint, a full SVD gives the basis $\mathbf Q_X$, so
$\operatorname{range}\mathbf Q_X=\ker\bm\Gamma_0^{\top}$
\cite{GolubVanLoan2013}.  When $\mathbf N_0$ is $\mathbf M$-orthonormal,
$\bm\Gamma_0$ also gives the coordinates of $\bm\Gamma$ in that basis;
otherwise it is a mass-coupling matrix with the same null-space condition.
Set $\mathbf Z_X=\mathbf N_0\mathbf Q_X$. Then
\[
  \operatorname{range}\mathbf Z_X
  =\left\{\mathbf N_0\bm y:
    \bm y\in\ker\bm\Gamma_0^{\top}\right\}.
\]
It follows that $\mathbf B_{\mathrm{bc}}\mathbf Z_X=0$ and
$\mathbf Z_X^{\top}\mathbf M\bm\Gamma=0$. If mass-inner-product
orthonormalization is applied,
$\mathbf Z_X^{\top}\mathbf M\mathbf Z_X=\mathbf I$; otherwise the reduced
mass matrix retains its general positive-definite form.

\subsection{Reduced Maxwell eigenproblem and workflow}

We discretize the weak Maxwell eigenproblem
\eqref{eq:continuous-weak-eigenproblem} on the physical space
\eqref{eq:abstract-physical-space} using the forms defined above. Write
$\bm c=\mathbf Z_X\bm y$ and test the sampled weak residual with
$\bm d=\mathbf Z_X\bm z$. Requiring it to vanish gives
\[
  \bm z^{\top}\mathbf Z_X^{\top}
  (\mathbf K-\lambda\mathbf M)\mathbf Z_X\bm y=0
  \quad\text{for every }\bm z.
\]
The resulting projected matrices have the Gram representations
\[
  \begin{aligned}
    \mathbf K_X
    &=\mathbf Z_X^{\top}\mathbf K\mathbf Z_X
      =(\widehat{\mathbf C}\mathbf Z_X)^{\top}
        (\widehat{\mathbf C}\mathbf Z_X),\\
    \mathbf M_X
    &=\mathbf Z_X^{\top}\mathbf M\mathbf Z_X
      =(\widehat{\bm\Phi}\mathbf Z_X)^{\top}
        (\widehat{\bm\Phi}\mathbf Z_X).
  \end{aligned}
\]
dR-Eig-PIELM solves
\[
  \mathbf K_X\mathbf Y=\mathbf M_X\mathbf Y\bm\Lambda,
  \quad
  \mathbf Y^{\top}\mathbf M_X\mathbf Y=\mathbf I,
\]
and recovers the vector-field coefficients as $\mathbf Z_X\mathbf Y$. Unlike
the baseline mixed strong-form product $\bm\Phi^{\top}\bm\Psi$, both reduced
matrices are symmetric Gram matrices.

\begin{algorithm}[H]
  \caption{Complete dR-Eig-PIELM workflow}
  \label{alg:drelm-complete}
  \renewcommand{\algorithmicrequire}{\textbf{Input:}}
  \renewcommand{\algorithmicensure}{\textbf{Output:}}
  \begin{algorithmic}[1]
    \REQUIRE $\widehat{\bm\Phi}$, $\widehat{\mathbf C}$,
    $\mathbf B_{\mathrm{bc}}$, $\bm\Gamma_p$, $\mathbf R_p$, $\mathbf T_0$;
    thresholds $\tau_b,\tau_g$; mode count $k$; optional
    $\mathbf M$-orthonormalization and guard $\tau_q$

    \STATE $\bm\Gamma\gets\bm\Gamma_p\mathbf T_0$
    \STATE require $\mathbf R_p\bm\Gamma_p=0$,
    $\mathbf C\bm\Gamma=0$, and
    $\mathbf B_{\mathrm{bc}}\bm\Gamma=0$; otherwise stop

    \STATE
    $\mathbf B_{\mathrm{bc}}
    =\mathbf U_b\bm\Sigma_b\mathbf V_b^{\top}$,
    $r_b=\operatorname{rank}_{\tau_b}(\mathbf B_{\mathrm{bc}})$,
    $\mathbf N_0=\mathbf V_b(:,r_b+1:\mathrm{end})$
    \hfill \textit{PEC null space}

    \IF{$\mathbf M$-orthonormalization is selected}
      \STATE
      $(\widehat{\bm\Phi}\mathbf N_0)\mathbf P_0
      =\mathbf Q_0\mathbf R_0$
      \hfill \textit{pivoted sample QR}
      \STATE require the nontruncating rank guard $\tau_q$ to pass;
      otherwise stop
      \STATE
      $\mathbf N_0\gets\mathbf N_0\mathbf P_0\mathbf R_0^{-1}$
    \ENDIF

    \STATE
    $\bm\Gamma_0
    \gets(\widehat{\bm\Phi}\mathbf N_0)^{\top}
    \widehat{\bm\Phi}\bm\Gamma$
    \STATE
    $\bm\Gamma_0
    =\mathbf U_g\bm\Sigma_g\mathbf V_g^{\top}$
    \STATE
    $r_g=\operatorname{rank}_{\tau_g}(\bm\Gamma_0)$,
    $\mathbf Q_X=\mathbf U_g(:,r_g+1:\mathrm{end})$
    \hfill \textit{gradient null space}
    \STATE require $r_g=N_s$; otherwise stop

    \STATE
    $\mathbf Z_X\gets\mathbf N_0\mathbf Q_X$
    \hfill \textit{physical-space basis}

    \STATE
    $\mathbf M_X\gets
    (\widehat{\bm\Phi}\mathbf Z_X)^{\top}
    (\widehat{\bm\Phi}\mathbf Z_X)$,
    $\quad
    \mathbf K_X\gets
    (\widehat{\mathbf C}\mathbf Z_X)^{\top}
    (\widehat{\mathbf C}\mathbf Z_X)$

    \STATE
    $\mathbf K_X\mathbf Y
    =\mathbf M_X\mathbf Y\bm\Lambda$,
    $\quad
    \mathbf Y^{\top}\mathbf M_X\mathbf Y=\mathbf I$

    \STATE retain the $k$ smallest eigenpairs $(\lambda_i,\mathbf y_i)$

    \ENSURE
    $\displaystyle
    \lambda_i,
    \bm E_i
    =\sum_j(\mathbf Z_X\mathbf y_i)_j\bm\varphi_j,
    \quad i=1,\ldots,k$
  \end{algorithmic}
\end{algorithm}

\section{Analysis of the compatible method}
\label{sec:compatible-analysis}

\subsection{Compatibility of the Construction}

Section~\ref{subsec:compatible-operator-realization} constructs the
three-dimensional tensor-product sequence and its two-dimensional
grad--curl form. We analyze their PEC-constrained polynomial spaces on
rectangles and cuboids \cite{ArnoldAwanou2014Cubical,BuffaEtAl2011IGAForms}.
For $d\in\{2,3\}$ and
$p\geq2$, let \(\mathbb P_r(I)\) denote the univariate polynomials of degree
at most \(r\) on an interval \(I\), and set
\begin{equation*}
  \Omega_d=\prod_{\alpha=1}^d(0,L_\alpha),
  \qquad
  \mathbb Q_{r_1,\ldots,r_d}(\Omega_d)
  =\bigotimes_{\alpha=1}^d\mathbb P_{r_\alpha}(0,L_\alpha),
\end{equation*}
and define
\begin{equation}\label{eq:dimension-independent-spaces}
\begin{aligned}
  \mathbb S_{p,d}^0
    &=\mathbb Q_{p,\ldots,p}(\Omega_d)\cap H_0^1(\Omega_d),\\
  \mathbb V_{p,d}
    &=\prod_{j=1}^d
      \mathbb Q_{p-\delta_{1j},\ldots,p-\delta_{dj}}(\Omega_d),\\
  \mathbb V_{p,d}^0
    &=\mathbb V_{p,d}\cap H_0(\curl;\Omega_d),
  \qquad \mathbb G_{p,d}=\nabla\mathbb S_{p,d}^0.
\end{aligned}
\end{equation}
Here $\delta_{ij}$ is the Kronecker delta. For $d=2$, these definitions
recover $\mathbb S_{p,p}^0$ in \eqref{eq:homogeneous-scalar-space} and the
vector trial space in \eqref{eq:graded-grad-curl-sequence}.
For $d=3$, $\mathbb V_{p,3}$ is the electric-field space in
\eqref{eq:three-dimensional-first-form-space} restricted to $\Omega_3$.
\begin{lemma}[Compatibility of the Tensor-Product Polynomial Construction]
	\label{lem:tensor-polynomial-compatibility}
	
	For each $d\in\{2,3\}$ and $p\geq2$, 
	$\mathbb G_{p,d}\subseteq\mathbb V_{p,d}^0$ and
	$(\nabla\times)\mathbb G_{p,d}=\{0\}$.
	Let $\gamma_\tau$ denote the tangential trace operator and define the
	tangential trace space
	$\mathcal T_{p,d}:=\{\gamma_\tau\mathbf v:
	\mathbf v\in\mathbb V_{p,d}\}$.
	Let $\mathcal S_{\partial}:\mathcal T_{p,d}\to\mathbb R^{N_b}$
	denote the boundary sampling operator. If $\mathcal S_{\partial}$ is
	injective on $\mathcal T_{p,d}$, then
	$\{\mathbf v_{\mathbf c}:\mathbf c\in\ker\mathbf B_{\mathrm{bc}}\}
	=\mathbb V_{p,d}^0$.
\end{lemma}

\begin{proof}
	Let $s\in\mathbb S_{p,d}^0$. Since
	$s\in\mathbb Q_{p,\ldots,p}(\Omega_d)$, differentiation with respect
	to $x_j$ gives
	$\partial_{x_j}s\in
	\mathbb Q_{p-\delta_{1j},\ldots,p-\delta_{dj}}(\Omega_d)$
	for $j=1,\ldots,d$. Hence $\nabla s\in\mathbb V_{p,d}$.
	
	Let $F$ be any boundary face of $\Omega_d$. Since $s$ is polynomial
	and $s|_F\equiv0$, for every tangential direction $\boldsymbol\tau$
	on $F$ we have
	$\boldsymbol\tau\cdot\nabla s
	=\partial_{\boldsymbol\tau}(s|_F)=0$.
	Thus the tangential trace of $\nabla s$ vanishes on every boundary
	face, that is, $\gamma_\tau(\nabla s)=0$ on $\partial\Omega_d$.
	Therefore $\nabla s\in\mathbb V_{p,d}^0$, and hence
	$\mathbb G_{p,d}\subseteq\mathbb V_{p,d}^0$.
	
	Moreover, since $s$ is polynomial, its mixed derivatives commute, so
	$\nabla\times\nabla s=0$. Since
	$\mathbb G_{p,d}=\nabla\mathbb S_{p,d}^0$, it follows that
	$(\nabla\times)\mathbb G_{p,d}=\{0\}$.
	
	It remains to identify the coefficient null space of
	$\mathbf B_{\mathrm{bc}}$ with $\mathbb V_{p,d}^0$.
	Let $\mathbf c\in\ker\mathbf B_{\mathrm{bc}}$. By the definition of
	the boundary constraint matrix,
	$\mathbf B_{\mathrm{bc}}\mathbf c
	=\mathcal S_{\partial}
	(\gamma_\tau\mathbf v_{\mathbf c})$.
	Thus $\mathbf B_{\mathrm{bc}}\mathbf c=0$ implies
	$\mathcal S_{\partial}
	(\gamma_\tau\mathbf v_{\mathbf c})=0$.
	Since $\gamma_\tau\mathbf v_{\mathbf c}\in\mathcal T_{p,d}$ and
	$\mathcal S_{\partial}$ is injective on $\mathcal T_{p,d}$, we obtain
	$\gamma_\tau\mathbf v_{\mathbf c}=0$ on $\partial\Omega_d$.
	Because $\mathbf v_{\mathbf c}\in\mathbb V_{p,d}$, this yields
	$\mathbf v_{\mathbf c}\in\mathbb V_{p,d}^0$.
	
	Conversely, if $\mathbf v_{\mathbf c}\in\mathbb V_{p,d}^0$, then
	$\gamma_\tau\mathbf v_{\mathbf c}=0$ on $\partial\Omega_d$.
	Applying the boundary sampling operator gives
	$\mathbf B_{\mathrm{bc}}\mathbf c
	=\mathcal S_{\partial}
	(\gamma_\tau\mathbf v_{\mathbf c})=0$,
	so $\mathbf c\in\ker\mathbf B_{\mathrm{bc}}$.
	Therefore
	$\{\mathbf v_{\mathbf c}:\mathbf c\in\ker\mathbf B_{\mathrm{bc}}\}
	=\mathbb V_{p,d}^0$.
\end{proof}

\subsection{Exactness and Ranks of the Construction}

\begin{theorem}[Exactness of the Discrete Grad--Curl Complex at
$\mathbb V_{p,d}^0$]
\label{thm:finite-exactness}
For each $d\in\{2,3\}$ and fixed $p\geq2$, the spaces in
\eqref{eq:dimension-independent-spaces} satisfy
\begin{equation}\label{eq:kernel-exactness}
  \ker\!\left(\left.\nabla\times\right|_{\mathbb V_{p,d}^0}\right)
  =\mathbb G_{p,d}.
\end{equation}
\end{theorem}

\begin{proof}
	Lemma~\ref{lem:tensor-polynomial-compatibility} gives
	$\mathbb G_{p,d}\subset
	\ker(\nabla\times|_{\mathbb V_{p,d}^0})$.
	For the reverse inclusion, let
	$\mathbf v\in\mathbb V_{p,d}^0$ satisfy
	$\nabla\times\mathbf v=0$, and fix
	$\mathbf x^\ast\in\Omega_d$.
	Since $\Omega_d$ is convex, define
	\[
	s(\mathbf x)
	=
	\int_0^1
	\mathbf v\bigl(\mathbf x^\ast+t(\mathbf x-\mathbf x^\ast)\bigr)
	\cdot(\mathbf x-\mathbf x^\ast)\,dt .
	\]
	Using $\partial_i v_j=\partial_j v_i$, we obtain
	$\partial_i s(\mathbf x)
	=\int_0^1 \frac{d}{dt}
	\bigl[t\,v_i(\mathbf x^\ast+t(\mathbf x-\mathbf x^\ast))\bigr]\,dt
	=v_i(\mathbf x)$.
	Hence $\nabla s=\mathbf v$.
	
	For the $j$th summand in the integrand,
	$v_j\in
	\mathbb Q_{p-\delta_{1j},\ldots,p-\delta_{dj}}$.
	Multiplication by $x_j-x_j^\ast$ raises the degree in $x_j$
	from at most $p-1$ to at most $p$, while all other coordinate
	degrees remain at most $p$. Integration in $t$ does not alter these
	degree bounds, so $s\in\mathbb Q_{p,\ldots,p}(\Omega_d)$.
	
	On each boundary face $F$, the PEC condition gives
	$\mathbf v_\tau=0$. Since $\nabla s=\mathbf v$,
	$\nabla_\tau s=0$ on $F$, and hence $s|_F$ is constant.
	Continuity of $s$ and the intersections of adjacent faces imply that
	these constants coincide on the connected boundary
	$\partial\Omega_d$. Subtracting this common constant leaves
	$\nabla s$ unchanged and yields
	$s\in\mathbb Q_{p,\ldots,p}(\Omega_d)\cap H_0^1(\Omega_d)
	=\mathbb S_{p,d}^0$.
	Thus $\mathbf v=\nabla s\in\mathbb G_{p,d}$, proving
	\[
	\ker\!\left(\nabla\times|_{\mathbb V_{p,d}^0}\right)
	=\mathbb G_{p,d}.
	\]
\end{proof}

\begin{corollary}[Ranks and dimensions of the tensor-polynomial reduction]
\label{cor:tensor-polynomial-ranks}
For fixed $p\geq2$ and $d\in\{2,3\}$, assume that the PEC trace samples are
unisolvent on the trace space of $\mathbb V_{p,d}$,
$\mathbf B_{\mathrm{bc}}\bm\Gamma=0$, and $\mathbf M$ is positive definite
on $\ker\mathbf B_{\mathrm{bc}}$.
In exact arithmetic,
\begin{equation}\label{eq:algebraic-dimension-chain}
\begin{aligned}
  \operatorname{rank}\mathbf B_{\mathrm{bc}}
    &=d\,p\bigl((p+1)^{d-1}-(p-1)^{d-1}\bigr),\\
  \operatorname{rank}\bm\Gamma
    &=N_{\nabla}=(p-1)^d,\\
  N_{\mathrm{phys}}
    &=(p-1)^{d-1}\bigl((d-1)p+1\bigr).
\end{aligned}
\end{equation}
For $d=2$, these three quantities are $4p$, $(p-1)^2$, and $p^2-1$,
respectively.
\end{corollary}

\begin{proof}
	For $s\in\mathbb S_{p,d}^0$, the zero trace on
	$x_\alpha=0,L_\alpha$ implies that
	$x_\alpha(L_\alpha-x_\alpha)$ divides $s$ for every $\alpha$. Hence
	\[
	\mathbb S_{p,d}^0
	=
	\left(\prod_{\alpha=1}^d x_\alpha(L_\alpha-x_\alpha)\right)
	\mathbb Q_{p-2,\ldots,p-2},
	\]
	and therefore $\dim\mathbb S_{p,d}^0=(p-1)^d$.
	
	For $\bm v\in\mathbb V_{p,d}^0$, the PEC condition implies that its
	$j$th component vanishes on the faces $x_\alpha=0,L_\alpha$ for
	$\alpha\ne j$. Thus
	$v_j$ contains the factor
	$\prod_{\alpha\ne j}x_\alpha(L_\alpha-x_\alpha)$, and hence
	$\dim(\mathbb V_{p,d}^0)_j=p(p-1)^{d-1}$. Consequently, the PEC and
	raw vector-space dimensions are
	\[
	\begin{aligned}
	N_{\mathrm{PEC}}&=\dim\mathbb V_{p,d}^0=d\,p(p-1)^{d-1},\\
	N_v&=\dim\mathbb V_{p,d}=d\,p(p+1)^{d-1}.
	\end{aligned}
	\]
	
	The gradient is injective on $\mathbb S_{p,d}^0$, so
	\[
	r_g:=\operatorname{rank}\bm\Gamma=\dim\mathbb G_{p,d}=(p-1)^d.
	\]
	Trace unisolvence gives
	$\ker\mathbf B_{\mathrm{bc}}\simeq\mathbb V_{p,d}^0$, and hence
	\[
	\operatorname{rank}\mathbf B_{\mathrm{bc}}=N_v-N_{\mathrm{PEC}}.
	\]
	The physical space is the mass-orthogonal complement of the represented
	gradient space within the PEC space. Positive definiteness of the mass
	inner product therefore gives
	\[
	N_{\mathrm{phys}}=N_{\mathrm{PEC}}-r_g.
	\]
	Substituting these counts yields \eqref{eq:algebraic-dimension-chain}.
	Setting $d=2$ gives the stated
	two-dimensional formulas.
\end{proof}

The rank and dimension formulas in Corollary~\ref{cor:tensor-polynomial-ranks}
hold in exact arithmetic. Numerical ranks coincide with the stated ranks when
the structural diagnostics detect no lost boundary or gradient directions.

\subsection{Well-posedness of the reduced problem}

\begin{corollary}[Well-posedness of the tensor-polynomial reduction]
\label{cor:tensor-polynomial-well-posedness}
For each $d\in\{2,3\}$ and fixed $p\geq2$, suppose that the PEC trace samples
are unisolvent, field evaluation is injective on $\mathbb V_{p,d}$, curl evaluation is
injective on the represented curl space, and the quadrature and material
weights are positive.  Then
$\mathbf K_X$ is symmetric positive definite.
Consequently, the reduced generalized eigenproblem is well posed and has a
positive real spectrum.
\end{corollary}

\begin{proof}
Trace unisolvence identifies the trial-function image of
$\ker\mathbf B_{\mathrm{bc}}$ with $\mathbb V_{p,d}^0$.  Field injectivity and positivity of
the quadrature and material weights give
\[
\mathbf M=\widehat{\bm\Phi}^{\mathsf T}\widehat{\bm\Phi}\succ0.
\]
By
Lemma~\ref{lem:tensor-polynomial-compatibility}, the trial-function image of
$\operatorname{range}\bm\Gamma$ is $\mathbb G_{p,d}$, with
\[
\mathbb G_{p,d}\subset\mathbb V_{p,d}^0.
\]
The null-space construction in Section~\ref{sec:algebraic-algorithm} then
gives a full-rank basis $\mathbf Z_X$ of $\mathbb X$, so
\begin{equation*}
  \mathbf M_X=\mathbf Z_X^{\mathsf T}\mathbf M\mathbf Z_X
  \succ0.
\end{equation*}

Let $\bm y$ be a reduced coefficient vector and let
$\bm v_{\bm y}\in\mathbb X$ be the trial field with full coefficient vector
$\mathbf Z_X\bm y$.  If $\bm y^{\mathsf T}\mathbf K_X\bm y=0$, the Gram
representation gives
\begin{equation*}
  0=\|\widehat{\mathbf C}\mathbf Z_X\bm y\|_2^2.
\end{equation*}
Positivity of the curl weights and injectivity of curl evaluation on the
represented curl space imply
\[
\nabla\times\bm v_{\bm y}=0.
\]
Because $\bm v_{\bm y}\in\mathbb X\subset\mathbb V_{p,d}^0$,
Theorem~\ref{thm:finite-exactness} gives $\bm v_{\bm y}\in\mathbb G_{p,d}$.
But $\mathbb X\subset\mathbb G_{p,d}^{\perp_m}$, and therefore
\[
  \bm y^{\mathsf T}\mathbf M_X\bm y
  =m(\bm v_{\bm y},\bm v_{\bm y})=0.
\]
Since $\mathbf M_X\succ0$, this implies $\bm y=0$, so $\mathbf K_X$ is
symmetric positive definite.  The finite-dimensional spectral theorem
gives a positive real spectrum.
\end{proof}

For a tensor Gauss grid with at least $p+1$ distinct nodes per coordinate,
successive univariate polynomial unisolvence makes both field evaluation and
curl evaluation injective on their stated represented spaces. For boundary
sampling, use at least $p+1$ distinct points on each edge in two dimensions.
In three dimensions, use a complete Cartesian product grid on every face,
with at least $p+1$ distinct nodes in each tangential coordinate, and sample
both tangential components. These boundary grids make trace sampling
injective and, together with the volume grid, satisfy the hypotheses above.
For a nonpolynomial material coefficient, positivity still establishes the
algebraic result, but the quadrature need not integrate the continuous
bilinear forms exactly.

\section{Numerical experiments}\label{sec:numerical-experiments}

The numerical experiments in this section test the structural identities of the compatible reduction developed in Sections~2--4, compare spectral accuracy and solver cost, and use mapped cutoff sensitivities in a local waveguide design.

All computations were performed on a 64-bit Windows~11 machine (build 26200) with an AMD Ryzen~9 8945HX processor and 32~GB of memory, using Python~3.12.12, NumPy~2.4.2, SciPy~1.16.3, and NGSolve~6.2.2606 \cite{HarrisEtAl2020NumPy,VirtanenEtAl2020SciPy,Schoeberl2014NGSolve}. OpenMP, BLAS, and NGSolve were restricted to a single thread before importing the numerical libraries. Timings reported for the rectangular cavity and annular-sector waveguide are averages over ten runs following one untimed warm-up. 

\subsection{Smooth asymmetric-material PEC cavity}
\label{sec:case1-structural-validation}

Let $\Omega=(0,1)^2$, $\epsilon=1$, and impose the PEC trace on the complete
boundary.  The inverse permeability is
\begin{equation*}
  \nu(x,y)=1+0.35\exp\!\left[-\left(\frac{x-0.63}{0.30}\right)^2
  -\left(\frac{y-0.41}{0.24}\right)^2\right].
\end{equation*}
Its off-center profile removes separability and the square symmetries, so no
analytic spectrum is available.

dR-Eig-PIELM uses the compatible space in \eqref{eq:graded-grad-curl-sequence} at $p=12$, with assembly by a $16\times16$ tensor-Gauss rule. The Eig-PIELM baseline uses degree-$12$ tensor-product Bernstein polynomials for both field components. Its trial space is $\mathbb S_{12,12}\times\mathbb S_{12,12}$. The baseline uses 2304 interior strong-form collocation points and 52 edgewise boundary samples (13 per edge), for 2356 samples in total. A second-order N\'ed\'elec method \cite{Nedelec1980,Hiptmair2002} supplies the reference spectrum on the finest 30,720-element mesh of a five-level refinement sequence.

\begin{figure}[htbp]
  \centering
  \includegraphics[width=\standardfigurewidth]{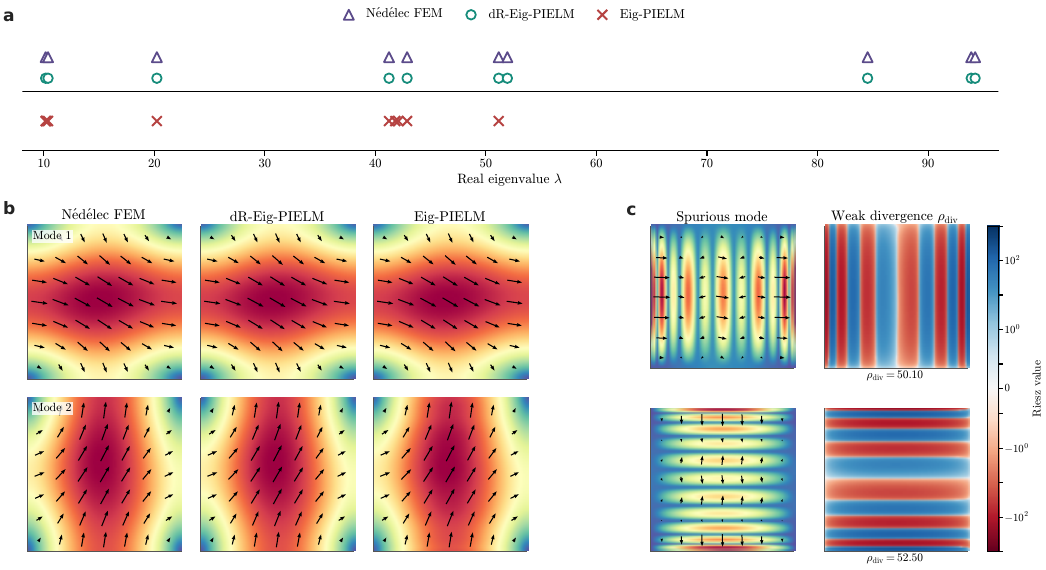}
  \caption{Spectral distribution and mode shapes of the smooth asymmetric-material PEC cavity. (a) First ten returned values.  (b) First two physical modes from N\'ed\'elec FEM, dR-Eig-PIELM, and Eig-PIELM. (c) Eig-PIELM fields at returned positions~2--3 and their weak-divergence Riesz representatives $r_{\rm div}$, with
  	$\rho_{\rm div}=\|r_{\rm div}\|_{L^2(\Omega)}/\|\bm E\|_{L^2(\Omega)}$.}
  \label{fig:case1-structural-summary-refactored}
\end{figure}
\FloatBarrier

Figure~\ref{fig:case1-structural-summary-refactored} compares the three methods'
spectra and representative fields.
The two spurious Eig-PIELM modes shown in
Figure~\ref{fig:case1-structural-summary-refactored}(c) have significant
gradient components and weak-divergence ratios of $50.10$ and $52.50$,
consistent with their exclusion from the discrete physical spectrum.
dR-Eig-PIELM removes the represented gradient directions through the
mass-orthogonal physical-space construction before the eigensolve.
Table~\ref{tab:case1-spectrum-refactored} lists the values underlying panel (a). dR-Eig-PIELM agrees with the first ten physical reference eigenvalues to amaximum relative difference of $6.713\times10^{-7}$. The baseline instead interleaves four spurious modes at positions 2, 3, 7, and 8.

\begin{table}[!htb]
	\caption{First ten physical reference and dR-Eig-PIELM eigenvalues,
		alongside the first ten values returned by Eig-PIELM. The symbol
		$^{\dagger}$ marks an Eig-PIELM mode classified as spurious because of
		its significant gradient component.}
	\label{tab:case1-spectrum-refactored}
	\centering
	\footnotesize
		\begin{tabular}{@{}crrr@{}}
			\toprule
			Returned position & FEM reference & dR-Eig-PIELM
			& Eig-PIELM \\
			\midrule
			1  & 10.165144 & 10.165144 & 10.165155 \\
			2  & 10.381589 & 10.381589 & $10.191586^{\dagger}$ \\
			3  & 20.227464 & 20.227464 & $10.272671^{\dagger}$ \\
			4  & 41.233365 & 41.233367 & 10.381638 \\
			5  & 42.867510 & 42.867513 & 20.227461 \\
			6  & 51.157612 & 51.157612 & 41.233548 \\
			7  & 51.940485 & 51.940488 & $41.800731^{\dagger}$ \\
			8  & 84.531219 & 84.531231 & $42.072129^{\dagger}$ \\
			9  & 93.908454 & 93.908517 & 42.867820 \\
			10 & 94.261127 & 94.261187 & 51.157539 \\
			\midrule
			\multicolumn{2}{r}{Maximum relative DR/FEM difference}
			& $6.713\times10^{-7}$ & \\
		
			\bottomrule
		\end{tabular}
\end{table}
\FloatBarrier

Table~\ref{tab:case1-rank-dimension-refactored} compares the computed ranks
and dimensions with Corollary~\ref{cor:tensor-polynomial-ranks} at $p=12$.
The first row counts independent tangential boundary constraints. The second
counts represented gradient directions that satisfy PEC, and the third gives
the dimension passed to the physical eigensolve.

\begin{table}[!htb]
\caption{Ranks and dimensions of the compatible reduction for the smooth
asymmetric-material PEC cavity at $p=12$.}
\label{tab:case1-rank-dimension-refactored}
\centering
\footnotesize
\begin{tabular}{@{}llr@{}}
\toprule
Quantity & Theoretical formula & Computed value \\
\midrule
PEC constraint rank & $\operatorname{rank}\mathbf B_{\mathrm{bc}}=4p$ & $48$ \\
Gradient-space dimension & $N_{\nabla}=(p-1)^2$ & $121$ \\
Physical-space dimension & $N_{\rm phys}=p^2-1$ & $143$ \\
\bottomrule
\end{tabular}
\end{table}
\FloatBarrier

\Needspace{6\baselineskip}
The PEC restriction reduces $312$ vector coefficients to $264$; removing
$121$ represented gradient directions leaves a $143$-dimensional physical
problem, as predicted by the corollary.

Table~\ref{tab:case1-structural-defects-refactored} reports algebraic residuals
for the represented gradient image, the physical basis, and the reduced
eigenproblem. Each of the four zero-product defects is measured by
$\|\mathbf A\mathbf B\|_F/(\|\mathbf A\|_F\|\mathbf B\|_F)$.
The last row reports
$\|\mathbf K_X\mathbf Y-\mathbf M_X\mathbf Y\bm\Lambda\|_F/
\|\mathbf K_X\|_F$.

The first two rows check the sampled PEC trace and curl-free property of the
represented gradients. The next two check that $\mathbf Z_X$ satisfies the
sampled PEC trace and is $\mathbf M$-orthogonal to the gradient image. The
reduced eigenpair residual is $1.342\times10^{-15}$. These residuals confirm
the tested finite-dimensional identities and algebraic solve at
floating-point roundoff level.

\begin{table}[htbp]
	\caption{Algebraic residuals of the dR-Eig-PIELM reduction for the smooth
		asymmetric-material PEC cavity.}
	\label{tab:case1-structural-defects-refactored}
	\centering
	\footnotesize
	\begin{tabular}{@{}llr@{}}
			\toprule
			Identity & Property & Result \\
			\midrule
			$\mathbf B_{\mathrm{bc}}\bm\Gamma=0$
			& PEC containment of represented gradients
			& $2.746\times10^{-17}$ \\
			
			$\mathbf C\bm\Gamma=0$
			& Curl-free gradient image
			& $3.106\times10^{-17}$ \\
			\midrule
			
			$\mathbf B_{\mathrm{bc}}\mathbf Z_X=0$
			& PEC admissibility of the physical basis
			& $4.031\times10^{-17}$ \\
			
			$\mathbf Z_X^{\top}\mathbf M\bm\Gamma=0$
			& Mass orthogonality to represented gradients
			& $6.540\times10^{-17}$ \\
			\midrule
			
			$\mathbf K_X\mathbf Y-\mathbf M_X\mathbf Y\bm\Lambda=0$
			& Reduced eigenpair equation
			& $1.342\times10^{-15}$ \\
			\bottomrule
	\end{tabular}
\end{table}
\FloatBarrier

To assess the effect of gradient removal, we also solve on the compatible PEC space with the represented gradient image retained. Both configurations use the same polynomial basis and PEC restriction. All $121=(p-1)^2$ represented gradient directions then form the numerical zero cluster. Removing them reduces the solved dimension from $264$ to $143$. The first ten positive eigenvalues change by at most $1.013\times10^{-13}$ relatively. The physical-space restriction removes the zero cluster before the eigensolve while preserving the positive spectrum to the reported precision. Although the mean eigensolve time decreases from $0.00536$ to $0.00169$~s, the end-to-end time increases from $0.01662$ to $0.03168$~s because, in this small-scale test, constructing the mass-orthogonal physical space costs more than is saved in matrix assembly and eigensolution.

\begin{samepage}
\subsection{Rectangular PEC cavity}
\label{sec:case2}

We consider the Maxwell eigenproblem~(1.1) on the rectangular cavity $\Omega=(0,1)\times(0,1.3)$ with $\epsilon=\nu=1$ and PEC walls. Separation of variables gives the positive eigenvalues \begin{equation*}
  \lambda_{mn}=\pi^2\left[m^2+\left(\frac{n}{1.3}\right)^2\right],
  \qquad m,n\geq0,\quad(m,n)\neq(0,0),
\end{equation*} where $m$ and $n$ index the $x$ and $y$ factors, respectively \cite{Collin1991,Monk2003}. Since $\lambda=\omega^2$ in~(1.1), these are squared resonant angular frequencies.
\end{samepage}

The finite-element method uses second-order triangular N\'ed\'elec $H(\curl)$ elements on four independently generated meshes. The tested dR-Eig-PIELM degrees are $p=5,\ldots,22$. Let $\mathcal I_{20}$ contain the twenty pairs $(m,n)$ with the smallest positive $\lambda_{mn}$. We compare each numerical eigenvalue in ascending order with the corresponding analytic label and measure
\begin{equation*}
  E_{\max}^{(20)}=
  \max_{(m,n)\in\mathcal I_{20}}
  \frac{|\lambda_{mn}^{\rm num}-\lambda_{mn}|}{\lambda_{mn}}.
\end{equation*}
Here $\lambda_{mn}^{\rm num}$ is paired with $\lambda_{mn}$ by spectral order. DOF counts the coefficients passed to the eigensolver. For FEM these are the free coefficients in the PEC N\'ed\'elec space; for dR-Eig-PIELM they are the reduced physical-space coefficients. We write $T_{\rm end}$ for the mean end-to-end solver time and set four targets
$\tau\in\{10^{-2},10^{-3},10^{-4},10^{-5}\}$ and compare configurations with $E_{\max}^{(20)}\leq\tau$. For each target, the table gives a second-order finite-element mesh and the lowest tested dR-Eig-PIELM degree that satisfy it. Each $T_{\rm end}$ measures the listed configuration alone and includes its mesh or basis construction, assembly, constraints, and eigensolve. Appendix~\ref{app:reproducibility-parameters} gives the remaining settings.

\begin{table}[!htbp]
\caption{Selected N\'ed\'elec FEM and dR-Eig-PIELM configurations at
four accuracy targets in the rectangular PEC cavity. Times are mean end-to-end
solver times in seconds; speedup is $T_{\rm FEM}/T_{\rm dR}$.}
\label{tab:case2-target-summary-refactored}
\centering
\footnotesize
\setlength{\tabcolsep}{4pt}
\begin{tabular}{@{}crrrrcrrrr@{}}
\toprule
& \multicolumn{4}{c}{N\'ed\'elec FEM ($r_{\rm FE}=2$)}
& \multicolumn{4}{c}{dR-Eig-PIELM} & \\
\cmidrule(lr){2-5}\cmidrule(lr){6-9}
Target $\tau$ & Triangles & DOF & $E_{\max}^{(20)}$ & $T_{\rm FEM}$
& $p$ & DOF & $E_{\max}^{(20)}$ & $T_{\rm dR}$ & Speedup \\
\midrule
$10^{-2}$ & 128 & 912 & $5.71\times10^{-3}$ & 0.208
& 10 & 99 & $1.23\times10^{-3}$ & 0.0133 & $15.6\times$ \\
$10^{-3}$ & 472 & 3453 & $4.89\times10^{-4}$ & 0.623
& 12 & 143 & $1.98\times10^{-5}$ & 0.0322 & $19.4\times$ \\
$10^{-4}$ & 1202 & 8877 & $7.33\times10^{-5}$ & 2.332
& 12 & 143 & $1.98\times10^{-5}$ & 0.0322 & $72.5\times$ \\
$10^{-5}$ & 3856 & 28674 & $7.32\times10^{-6}$ & 11.471
& 14 & 195 & $1.51\times10^{-7}$ & 0.0505 & $227.4\times$ \\
\bottomrule
\end{tabular}
\end{table}

Table~\ref{tab:case2-target-summary-refactored} compares the selected finite-element and dR-Eig-PIELM configurations at each target. For these configurations, the end-to-end speedup rises from $15.6\times$ at $10^{-2}$ to $227.4\times$ at $10^{-5}$. The $p=12$ dR-Eig-PIELM configuration satisfies both the $10^{-3}$ and $10^{-4}$ targets, so its measured result appears in both rows.

With the same global polynomial basis, we test quadrature order and numerical
rank tolerances at $p=10,12,14,16$.
First, we fix $\tau_b=10^{-12}$ and $\tau_g=10^{-11}$ and vary
$q-p\in\{1,2,3,4,5\}$. Next, we fix $q=p+4$ and vary
$\tau_b=s_\tau 10^{-12}$ and $\tau_g=s_\tau 10^{-11}$, with
$s_\tau\in\{10^{-2},10^{-1},1,10,10^2\}$.

\begin{figure}[htbp]
  \centering
  \includegraphics[width=\standardfigurewidth]{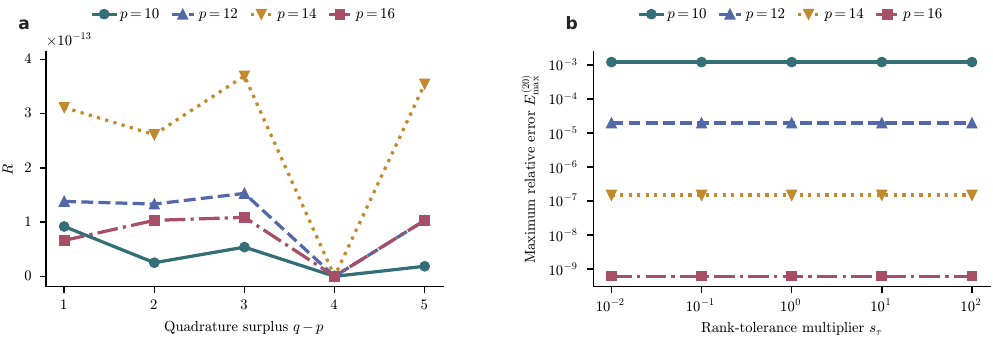}
  \caption{Case~2 parameter sensitivity for dR-Eig-PIELM at
  $p=10,12,14,16$: (a) quadrature order margin $q-p$ at fixed $s_\tau=1$;
  (b) rank-tolerance multiplier $s_\tau$ at fixed $q=p+4$.
  Panel (a) shows the maximum spectral drift $R$ of the first twenty positive
  eigenvalues relative to $(q-p,s_\tau)=(4,1)$.
  Panel (b) shows their maximum relative error $E_{\max}^{(20)}$.}
  \label{fig:case2-parameter-robustness}
\end{figure}

Figure~\ref{fig:case2-parameter-robustness} summarizes the quadrature and
rank-tolerance sensitivity results. At each degree,
$R=\max_{1\le j\le20}|\lambda_j^{\rm test}-
\lambda_j^{\rm base}|/\lambda_j^{\rm base}$.  The largest quadrature-scan
drift is $3.69\times10^{-13}$; the tolerance scan gives $R=0$ at the reported
precision.  All 36 distinct settings pass the checks and retain
$N_{\rm phys}=99,143,195,255$ for $p=10,12,14,16$, respectively.

We finally compare the numerical effect of post-PEC mass orthonormalization
with a smoothly varying material coefficient. This test
retains the Case~2 rectangle, the normalized shifted-Legendre basis, and the
same PEC construction, but replaces the unit permittivity used above by
$\epsilon(x,y)=10^{2(x-1)}$ with $\nu=1$, so that the permittivity ranges
from $0.01$ to $1$.
Both arms use a $36\times36$ tensor-Gauss rule,
$\tau_b=10^{-12}$, and $\tau_g=10^{-11}$ for $p=5,\ldots,20$. The \emph{raw}
arm retains the PEC null-space coordinates, whereas the \emph{stable} arm
mass-orthonormalizes those same coordinates after the PEC restriction.

\begin{figure}[htbp]
	\centering
	\includegraphics[width=\standardfigurewidth]{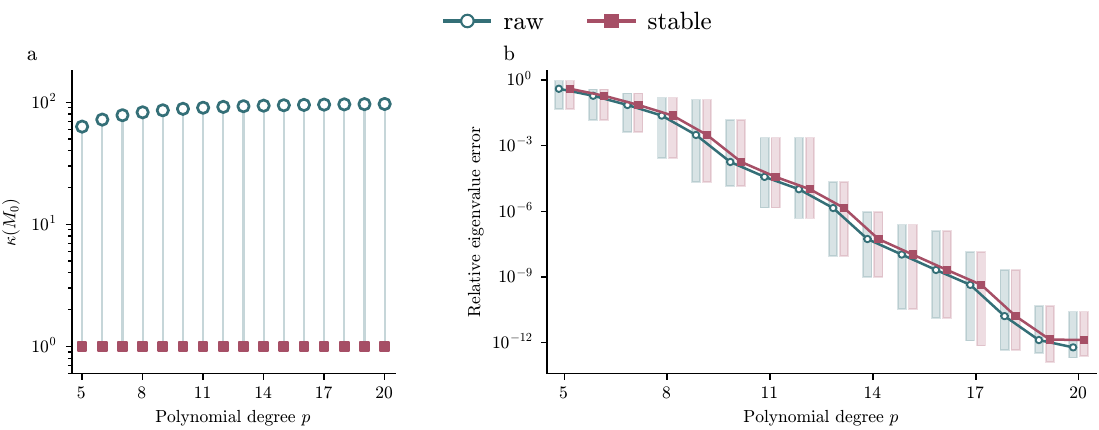}
	\caption{Post-PEC orthonormalization comparison in Case~2 with
		$\epsilon(x,y)=10^{2(x-1)}\in[0.01,1]$ and $\nu=1$.
		Raw and stable differ only in whether the PEC null-space coordinates are
		mass-orthonormalized after PEC restriction. (a) Condition number of the
		post-PEC weighted Gram matrix. (b) Relative errors of the first twenty
		positive eigenvalues. Boxes show the 25th--75th percentiles, representing
		the central 50\% of the twenty errors; lines connect the medians.}
	\label{fig:case2-post-pec-stress}
\end{figure}

Figure~\ref{fig:case2-post-pec-stress}(a) shows the condition number of the
post-PEC weighted Gram matrix. At $p=20$, its condition number is $96.85$
for raw and approximately one for stable. Both realizations retain the
expected gradient rank $(p-1)^2$ throughout the scan. The first twenty positive eigenvalues are
compared with an independent separated-variable reference evaluated at
85-digit precision. Figure~\ref{fig:case2-post-pec-stress}(b) shows the
central quartiles and medians of the twenty modewise relative errors,
with comparable degree convergence for the two realizations. At $p=20$,
both maximum relative
errors are $3.33\times10^{-6}$. Thus, post-PEC orthonormalization improves
the mass-coordinate conditioning while preserving comparable spectral
accuracy in this test.

\begingroup
\raggedbottom
\subsection{Annular-sector PEC waveguide}
\label{sec:case3}

We optimize the cross-sectional area of an air-filled PEC annular-sector
channel, a local model of a $3.6\,\mathrm{GHz}$ sectoral-waveguide mode
converter \cite{Kumar2017Sectoral}. Related all-metal converters have been
tested experimentally \cite{Kumar2019Sectoral}.

The initial cross-section is
\begin{equation*}
  \Omega_s=\{(r,\theta):r_i<r<r_o,\ 0<\theta<\alpha\},
  \quad r_i=0.030\ \mathrm m,\quad r_o=0.060\ \mathrm m,
  \quad \alpha=\pi/3.
\end{equation*}
With the free-space material constants factored out, the cross-sectional
matrices use normalized coefficients $\epsilon=\nu=1$. Here
$\lambda_j=k_{c,j}^2$ is the squared cutoff wavenumber, and the cutoff
frequency is
\[
  f_{c,j}=\frac{c_0}{2\pi}\sqrt{\lambda_j},
  \qquad c_0=299\,792\,458\ \mathrm{m\,s^{-1}}.
\]
We label separated cutoff modes as $\mathrm{TE}_{mn}$: $m\geq0$ is the
angular index in $\cos(m\pi\theta/\alpha)$, and $n\geq1$ counts the positive
radial Neumann--Bessel roots.
At $f_{\rm op}=3.6$~GHz, its initial cutoffs satisfy
\[
  f_{c,1}=f_{c,11}=3.1473\ \mathrm{GHz}
  <3.6\ \mathrm{GHz}
  <f_{c,2}=f_{c,01}=5.0840\ \mathrm{GHz}.
\]
The $\mathrm{TE}_{11}$ mode propagates and $\mathrm{TE}_{01}$ remains
evanescent. We freeze
this cutoff pair before optimization and require every candidate to retain it.
Figure~\ref{fig:case3-geometry-mode} shows the initial geometry and its
dominant transverse mode.
With $r_o=60$~mm fixed, the compactness problem is
\begin{equation}\label{eq:case3-compact-design}
  \begin{aligned}
  \min_{r_i,\alpha}\quad
    &A(r_i,\alpha)=\frac{\alpha}{2}\left(r_o^2-r_i^2\right),\\
  \text{subject to}\quad
    &f_{c,1}(r_i,\alpha)\le f_{c,1}^{0},\qquad
     f_{c,2}(r_i,\alpha)\ge f_{c,2}^{0},\\
    &0.95\le r_i/r_i^0\le1.05,\qquad
     0.95\le\alpha/\alpha^0\le1.05,
  \end{aligned}
\end{equation}
where $0$ denotes the initial design. The cutoff constraints preserve its
spectral envelope as the cross-sectional area decreases.

\begin{figure}[htbp]
  \centering
  \includegraphics[width=0.8\linewidth]{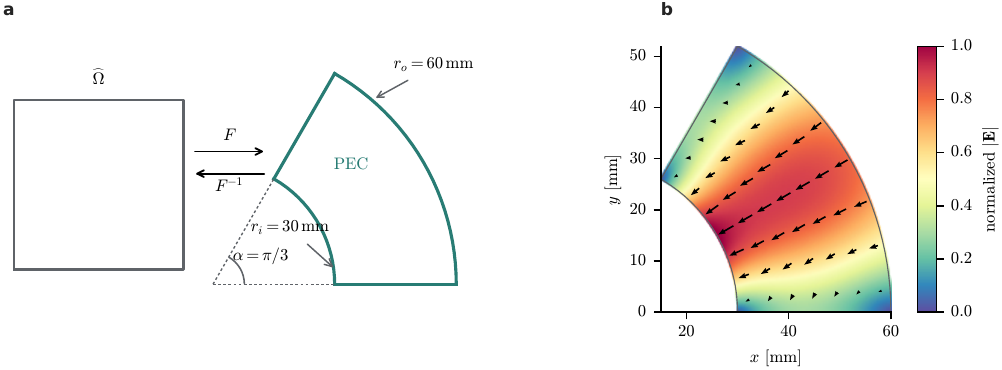}
  \caption{Nominal annular-sector PEC waveguide. (a) Mapping from the reference
  square $\widehat\Omega$ to the cross-section with $r_i=30$~mm,
  $r_o=60$~mm, and $\alpha=\pi/3$. (b) Normalized magnitude and transverse
  direction of the dR-Eig-PIELM $\mathrm{TE}_{11}$ electric-field mode.}
  \label{fig:case3-geometry-mode}
\end{figure}

At each trial geometry, the reference square $\widehat\Omega=(0,1)^2$ is
mapped to the annular sector by
\begin{equation*}
  r=r_i+(r_o-r_i)\xi,\qquad
  \theta=\alpha\eta,\qquad
  F(\xi,\eta)=(r\cos\theta,r\sin\theta).
\end{equation*}
With $\mathsf J=DF$, the vector field is transferred covariantly,
\begin{equation}\label{eq:case3-covariant-map-refactored}
  \bm E\circ F=\mathsf J^{-\top}\widehat{\bm E},
  \qquad
  (\nabla\times\bm E)\circ F
  =\frac{1}{\det \mathsf J}\widehat\nabla\times\widehat{\bm E}.
\end{equation}
This transfer preserves the PEC tangential trace and the represented
grad--curl relation \cite{RognesKirbyLogg2009,Monk2003}.

At each design iterate, dR-Eig-PIELM assembles the mapped matrices and
computes the two lowest positive cutoffs. Sensitivities are evaluated in
the geometry-independent reference PEC space, whose positive spectrum
agrees with that of the physical-space reduction. For
$\mu\in\{r_i,\alpha\}$, let $\bm u_j$ denote the full coefficient vector
of a simple positive eigenmode. Its sensitivity satisfies
\begin{equation*}
  \frac{\partial\lambda_j}{\partial\mu}
  =\bm u_j^{\top}\!\left(
  \frac{\partial\mathbf K}{\partial\mu}
  -\lambda_j\frac{\partial\mathbf M}{\partial\mu}\right)\!\bm u_j,
  \qquad \bm u_j^{\top}\mathbf M\bm u_j=1.
\end{equation*}
This is the standard first-order sensitivity formula for a simple symmetric
generalized eigenpair \cite{FoxKapoor1968}.
The corresponding cutoff derivative is
\[
  \frac{\partial f_{c,j}}{\partial\mu}
  =\frac{c_0}{4\pi\sqrt{\lambda_j}}
    \frac{\partial\lambda_j}{\partial\mu}.
\]
The mapped quadrature weights and the area objective are differentiated
analytically.  Modes are matched by mass-weighted overlap.  A $21\times21$
grid of complete dR-Eig-PIELM solves identifies a feasible starting
  point.  SLSQP, a sequential quadratic programming method
  \cite{NocedalWright2006}, then minimizes the area under the two cutoff
  constraints.

Against the analytic Bessel reference \cite{Collin1991}, the maximum
relative errors over the first ten positive cutoffs at the optimized geometry
are $3.777\times10^{-9}$ for dR-Eig-PIELM and $3.001\times10^{-10}$ for
the independent geometry-fitted N\'ed\'elec control. The full reduced solve
and the reference PEC eigenproblem agree within $7.02\times10^{-14}$~GHz
for the first two cutoffs. At the initial geometry, nonzero cutoff directional
derivatives agree with centered finite differences to a maximum relative
discrepancy of $1.24\times10^{-6}$. Appendix~\ref{app:reproducibility-parameters}
gives the discretization and verification parameters.

The optimized design, $r_i=31.5$~mm and $\alpha=59.243074^\circ$,
reduces the area from $1413.717$ to $1348.189$~mm$^2$ ($4.635\%$).
The $\mathrm{TE}_{11}$ constraint is active at $3.147301881$~GHz;
the $\mathrm{TE}_{01}$ cutoff rises to $5.339524155$~GHz.
Table~\ref{tab:case3-compact-design} summarizes the design and constraint checks.

\begin{table}[H]
\caption{Compactness optimization with $r_o=60$~mm. Cutoffs are verified by
complete dR-Eig-PIELM re-solves; parentheses give changes from the initial design.}
\label{tab:case3-compact-design}
\centering
\footnotesize
\begin{tabular}{@{}llll@{}} 
	\toprule 
	Quantity & Design requirement & Initial & Optimized \\ 
	\midrule 
	$r_i$ (mm) 
	& $28.5\le r_i\le31.5$ 
	& 30.000 
	& 31.500 \\ 
	
	$\alpha$ (deg) 
	& $57\le\alpha\le63$ 
	& 60.000 
	& 59.243 \\ 
	
	$A$ (mm$^2$) 
	& minimize 
	& 1413.717 
	& 1348.189 ($-65.527$) \\ 
	
	$f_{c,11}$ ($\mathrm{TE}_{11}$, GHz) 
	& $\le3.147302$ 
	& 3.147302 
	& 3.147302 \\ 
	
	$f_{c,01}$ ($\mathrm{TE}_{01}$, GHz) 
	& $\ge5.083993$ 
	& 5.083993 
	& 5.339524 ($+0.255531$) \\ 
	
	Operating condition 
	& $f_{c,11}<3.6<f_{c,01}$ 
	& satisfied 
	& satisfied \\ 
	\bottomrule 
\end{tabular}

\end{table}

The inner radius reaches its prescribed upper bound, and the active first
cutoff constraint determines the final angle. This local cross-sectional
design retains the frozen spectral specification while reducing area.
dR-Eig-PIELM supplies both cutoff values and their gradients within the
constrained geometry-design loop.
\FloatBarrier
\endgroup

\section{Conclusion}\label{sec:conclusion}

dR-Eig-PIELM extends the linear-coefficient and one-shot architecture of Eig-PIELM to PEC Maxwell eigenproblems by combining compatible tensor-product grad--curl spaces with a mass-orthogonal physical-space restriction. This construction removes the represented gradient image from the sampled PEC null space before the eigensolve. For the tensor-polynomial spaces considered here, fixed-degree exactness identifies that image with the discrete curl kernel. Under the stated sampling and material-positivity assumptions, the reduced eigenproblem is symmetric-definite.
In the numerical experiments, the method excludes the spurious eigenmodes with significant gradient components admitted by the boundary-only baseline. The smooth test problems show rapid degree convergence and favorable work--precision for the compared configurations. A covariantly mapped waveguide example retains spectral accuracy and uses validated cutoff sensitivities in a local design update.
Further work should develop a full spectral-convergence theory, including uniformly stable commuting projections and discrete-compactness estimates. Practical three-dimensional Maxwell problems also require replacing the present dense global constraint reductions with sparse, localized, or block-structured constructions.

\section*{Acknowledgments} The authors would like to express sincere gratitude for the valuable discussions and support from professor St\'ephane Lanteri, Inria.)

\section*{Declaration of generative AI and AI-assisted technologies in the manuscript preparation process}
During the preparation of this work, the authors used ChatGPT by OpenAI to
support language editing, improve readability, and obtain feedback on
manuscript structure and phrasing. The authors reviewed, edited, and verified
the content as needed and take full responsibility for the content of the
published article.

\newpage

\appendix
\FloatBarrier
\section{Reproducibility parameters}
\label{app:reproducibility-parameters}

This appendix lists only the implementation settings not stated in
Section~\ref{sec:numerical-experiments}.  Unless varied in the untimed
Case~2 audit, every two-dimensional dR-Eig-PIELM run uses relative SVD
thresholds $\tau_b=10^{-12}$ for the sampled PEC matrix and
$\tau_g=10^{-11}$ for the represented gradient coordinates. Case~1 uses the
mass-coupling null-space construction without post-PEC mass
orthonormalization, as do the Case~2 runs in
Table~\ref{tab:case2-target-summary-refactored}. The Case~2 parameter study
and all Case~3 runs apply post-PEC mass orthonormalization. In those runs,
the pivoted sample QR is
nontruncating, with rank guard
$\epsilon_{\rm mach}\max(N_{\rm row},N_{\rm col})$ for the sample-QR matrix
$\widehat{\bm\Phi}\mathbf N_0$.  Polynomial edge traces use
$p+1$ equispaced points per edge.  A failed guard terminates the configuration
rather than reducing the trial space.
The variable-material orthonormalization comparison in
Fig.~\ref{fig:case2-post-pec-stress} pairs the first twenty values with the independent
reference in ascending spectral order. 

\begin{table}[H]
\centering
\caption{Additional case-specific settings not stated in the numerical-experiment text.
Here $q$ is the one-dimensional Gauss order, $r_{\rm FE}$ is the N\'ed\'elec
order, and $\tau_+$ is the positive-eigenvalue threshold for solves retaining
the gradient kernel.}
\label{tab:case-parameters}
\footnotesize
\setlength{\tabcolsep}{4pt}
\begin{tabular}{>{\raggedright\arraybackslash}p{.07\linewidth}
>{\raggedright\arraybackslash}p{.39\linewidth}
>{\raggedright\arraybackslash}p{.47\linewidth}}
\toprule
Case & dR-Eig-PIELM / Eig-PIELM & N\'ed\'elec control \\
\midrule
1 & No post-PEC mass orthonormalization; Eig-PIELM uses $\tau_+=10^{-7}$ and PEC SVD
threshold $\tau_b=10^{-12}$ & NGSolve \texttt{maxh}=0.28, bonus integration order 8,
exact gradient projector, and fixed 80-step PINVIT \\
2 & Table~\ref{tab:case2-target-summary-refactored}: no post-PEC mass
orthonormalization; $q=p+4$ & NGSolve
\texttt{maxh}=0.15, 0.08, 0.05, 0.028 for the four targets, respectively;
$r_{\rm FE}=2$, meshes generated directly; bonus integration order 8;
exact gradient projector; residual tolerance
$10^{-4}$ checked every 5 steps, cap 600, seed 1729 \\
3 & $p=6,8,10,12,14$, $q=p+4$; fixed-PEC sensitivity solves use $\tau_+=10^{-5}$ & Curved
NGSolve with $r_{\rm FE}=1,\ldots,6$ at $h=8$~mm and $r_{\rm FE}=4$ at
$h=12,8,6,4$~mm; geometry order $\max(3,r_{\rm FE}+1)$; bonus integration
order $r_{\rm FE}+4$; exact gradient constraint; eigensolver tolerance
$10^{-10}$ \\
\bottomrule
\end{tabular}
\end{table}

The Case~1 quotient ablation uses one unrecorded warm-up for each arm and 20
paired repetitions.  The within-pair order alternates between
repetitions, and the reported times are sample means.  Timed regions exclude
diagnostics, reference-error evaluation, file output, and plotting.

The untimed Case~2 scans use $p=10,12,14,16$.  With
$\tau_b=10^{-12}$ and $\tau_g=10^{-11}$ fixed, it varies $q-p$ over
$1,2,3,4,5$.  At fixed $q=p+4$, the tolerance scan multiplies both thresholds
by $s_\tau=10^{-2},10^{-1},1,10,10^2$.  The shared baseline is evaluated
once, giving nine distinct settings per degree.

The Case~3 design runs use $p=14$ and $q=18$.  The terminal design is checked
by a complete dR-Eig-PIELM re-solve and a curved sixth-order N\'ed\'elec solve at
$h=8$~mm.  The derivative checks use relative steps
$10^{-2},\ldots,10^{-6}$ in the two coordinate directions and three fixed
random directions; modes are matched by mass-weighted overlap.

\section{Coefficient realization of compatible trial spaces}
\label{app:three-dimensional-coefficient-realization}

Let $\mathbf D_p^{(\alpha)}$ be the exact univariate derivative matrix from
degree $p$ to degree $p-1$ in direction $\alpha$. On an interval of length
$L_\alpha$, the physical derivative matrix is
$\widetilde{\mathbf D}_p^{(\alpha)}
=L_\alpha^{-1}\mathbf D_p^{(\alpha)}$ when $\mathbf D_p^{(\alpha)}$
is formed in the reference coordinate. Let $\mathbf I_r^{(\alpha)}$
denote the identity on the degree-$r$ coefficient space. For the Bernstein
basis, differentiation follows from
$\partial_\xi B_i^p=p(B_{i-1}^{p-1}-B_i^{p-1})$
\cite{Farouki2012Bernstein}; shifted-Legendre coordinates use exact
polynomial derivative expansions.

\subsection{Two-dimensional maps}

On $\Omega_2=(0,L_x)\times(0,L_y)$, the scalar space is
$\mathbb S_{r,s}=\mathbb P_r(0,L_x)\otimes\mathbb P_s(0,L_y)$.
With componentwise ordering of the vector coefficients, the gradient and
scalar-curl matrices for \eqref{eq:graded-grad-curl-sequence} are
\begin{equation*}
  \bm\Gamma_p=
  \begin{bmatrix}
    \widetilde{\mathbf D}_p^{(x)}\otimes\mathbf I_p^{(y)}\\
    \mathbf I_p^{(x)}\otimes\widetilde{\mathbf D}_p^{(y)}
  \end{bmatrix},
  \qquad
  \mathbf R_p=
  \begin{bmatrix}
    -\mathbf I_{p-1}^{(x)}\otimes\widetilde{\mathbf D}_p^{(y)}
    &\widetilde{\mathbf D}_p^{(x)}\otimes\mathbf I_{p-1}^{(y)}
  \end{bmatrix}.
\end{equation*}
The tensor factors give $\mathbf R_p\bm\Gamma_p=0$. If
$\mathbf N_x^0$ and $\mathbf N_y^0$ span the endpoint-null univariate
scalar spaces, the PEC scalar injection and restricted gradient are
\begin{equation*}
  \mathbf T_0=\mathbf N_x^0\otimes\mathbf N_y^0,
  \qquad \bm\Gamma=\bm\Gamma_p\mathbf T_0.
\end{equation*}

\subsection{Three-dimensional maps}
\label{subsec:three-dimensional-realization}

On $\Omega_3=(0,L_x)\times(0,L_y)\times(0,L_z)$, the unrestricted
gradient map of \eqref{eq:graded-three-dimensional-sequence} is
\begin{equation*}
  \bm\Gamma_{p,3}=
  \begin{bmatrix}
    \widetilde{\mathbf D}_p^{(x)}\otimes\mathbf I_p^{(y)}\otimes\mathbf I_p^{(z)}\\
    \mathbf I_p^{(x)}\otimes\widetilde{\mathbf D}_p^{(y)}\otimes\mathbf I_p^{(z)}\\
    \mathbf I_p^{(x)}\otimes\mathbf I_p^{(y)}\otimes\widetilde{\mathbf D}_p^{(z)}
  \end{bmatrix}.
\end{equation*}
The PEC scalar injection is
$\mathbf T_{0,3}=\mathbf N_x^0\otimes\mathbf N_y^0\otimes\mathbf N_z^0$,
and $\bm\Gamma_3=\bm\Gamma_{p,3}\mathbf T_{0,3}$ represents
$\mathbb G_{p,3}=\nabla\mathbb S_{p,3}^0$.
The curl matrix applies the signed directional derivatives
\begin{equation*}
  \mathbf R_{p,3}
  \begin{bmatrix}\bm v_x\\\bm v_y\\\bm v_z\end{bmatrix}
  =\begin{bmatrix}
    \mathbf D_y\bm v_z-\mathbf D_z\bm v_y\\
    \mathbf D_z\bm v_x-\mathbf D_x\bm v_z\\
    \mathbf D_x\bm v_y-\mathbf D_y\bm v_x
  \end{bmatrix},
\end{equation*}
where each $\mathbf D_\alpha$ acts on the corresponding tensor factor
with the physical interval scaling. Consequently,
$\mathbf R_{p,3}\bm\Gamma_{p,3}=0$. The PEC matrix samples both
tangential components on every face on a complete Cartesian product grid,
with at least $p+1$ distinct nodes in each tangential coordinate.

Let $\mathbf W_3$ be the tensor Gauss-weight matrix, let $\mathbf I_3$ be the
$3\times3$ identity on field components, and let
$\mathbf N_\epsilon,\mathbf N_\nu$ contain material values at the quadrature
nodes.  Denote by $\bm\Phi_3$ and $\mathbf C_3$ the field- and vector-curl
evaluation matrices, with their components stacked consistently.
Define
\begin{equation*}
  \mathbf W_{\epsilon,3}
  =\mathbf I_3\otimes(\mathbf W_3\mathbf N_\epsilon),
  \qquad
  \mathbf W_{\nu,3}
  =\mathbf I_3\otimes(\mathbf W_3\mathbf N_\nu),
\end{equation*}
and
\begin{equation*}
  \widehat{\bm\Phi}_3=\mathbf W_{\epsilon,3}^{1/2}\bm\Phi_3,
  \qquad
  \widehat{\mathbf C}_3=\mathbf W_{\nu,3}^{1/2}\mathbf C_3.
\end{equation*}
Together with the three-dimensional $\mathbf B_{\mathrm{bc}}$ and
$\bm\Gamma_3$, these weighted evaluation matrices provide the field,
curl, boundary, and gradient inputs to Algorithm~\ref{alg:drelm-complete}.

\bibliographystyle{siamplain}
\bibliography{references}

\end{document}